\documentclass{amsart}
\usepackage{amssymb}
\usepackage{amsmath}
\usepackage[centertags]{amsmath}
\usepackage[english,francais]{babel}

\newtheorem{theorem}{Theorem}[section]
\newtheorem{lemma}[theorem]{Lemma}
\newtheorem{proposition}[theorem]{Proposition}
\newtheorem{corollary}[theorem]{Corollary}
\theoremstyle{definition}
\newtheorem{definition}[theorem]{Definition}
\newtheorem{example}[theorem]{Examples}

\theoremstyle{remark}
\newtheorem{remark}[theorem]{Remark}

\numberwithin{equation}{section}

\begin{document}

\title[ Extremal almost-K\"ahler metrics  ]{\bf  Extremal almost-K\"ahler metrics  }
\author{Mehdi Lejmi}

\address{ D{\'e}partement de Math{\'e}matiques\\
UQAM\\ C.P. 8888 \\ Succ. Centre-ville \\ Montr{\'e}al (Qu{\'e}bec) \\
H3C 3P8 \\ Canada} \email{lejmi.mehdi@courrier.uqam.ca}
%\thanks{The author thanks Prof. V. Apostolov for his invaluable help and judicious advice.}

\maketitle
\medskip

\selectlanguage{english}
\begin{abstract}
We generalize the notions of the Futaki invariant and extremal vector field of a compact K\"ahler manifold to the general almost-K\"ahler case and show the periodicity of the extremal vector field when the symplectic form represents an integral cohomology class modulo torsion. We also give an explicit formula of the hermitian scalar curvature in Darboux coordinates which allows us to obtain examples of non-integrable extremal almost-K\"ahler metrics saturating LeBrun's estimates.
\end{abstract}

\section{Introduction}
Let $(M^{2n},\omega)$ be a compact symplectic manifold of dimension $2n$. Recall that an almost-complex structure $J$ is compatible with $\omega$ if the tensor field $g(\cdot,\cdot)= \omega(\cdot,J\cdot)$ defines a Riemannian metric on
$M$; in this case, the triple $(\omega,J,g)$ is referred to as an ($\omega$-compatible) almost-K\"ahler structure on $M$. Any such
structure defines, in a canonical way, a hermitian connection $\nabla$ on the complex tangent bundle $(T(M),J,g)$.
%which induces the usual 
%Cauchy-Remann operator $\bar \partial$ on $T(M)$, i.e. such that $ 
%\nabla_{X^{0,1}}Y^{1,0} = [X^{0,1}, Y^{1,0}]^{1,0}$ for any complex 
%vector fields $X^{0,1}$ and $Y^{1,0}$ of types $(0,1}$ and $(1,0)$, 
%respectively. Note that $\nabla$ has torsion unless $J$ is integrable, 
%i.e. when $(g,J,\omega)$ is a K\"ahler structure. 
Taking trace and contracting the curvature of $\nabla$ by $\omega$, one obtains the
hermitian scalar curvature $s^{\nabla}$ of $(\omega,J,g)$.

It is well-known \cite{don,fuj} that the space of all $\omega$-compatible almost-K\"ahler structures, here denoted by 
$AK_{\omega}$, is a contractible Fr\'echet manifold endowed with a formal K\"ahler structure. The infinite dimensional group $\textrm{Ham}(M,\omega)$ of hamiltonian symplectomorphisms naturally acts on
$AK_{\omega}$ and a crucial observation of Donaldson~\cite{don} (generalizing \cite{fuj} to the non-integrable almost-K\"ahler case)
is that this action is hamiltonian with moment map $\mu : AK_{\omega} \to (\textrm{Lie}(\textrm{Ham}(M,\omega))^\ast$ given by
$\mu_J(f) = \int_M s^{\nabla}f \frac{\omega^n}{n!},$
where $f$ is any smooth function with zero
integral on $M$, viewed also as an element of the Lie algebra of
$\textrm{Ham}(M,\omega)$. As already observed in \cite{apo-dra}, this interpretation
of $s^{\nabla}$ immediately implies that the critical points of the
functional $J \mapsto \int_M \left({s^{\nabla}}\right)^2 \frac{\omega^n}{n!}$ over $AK_{\omega}$ are the almost-K\"ahler metrics for which the symplectic gradient
of the hermitian scalar curvature is an infinitesimal isometry of the almost-complex structure $J$.
%satisfying the Euler-Lagrange equation $\mathfrak{L}_{grad_{\omega} s^{\nabla}} J =0.$
Motivated by striking analogy with the notion of {\it extremal K\"ahler metric} introduced by Calabi~\cite{cal}, we refer to the almost-K\"ahler metrics
verifying the above condition as {\it extremal almost-K\"ahler metrics.}

The above formal picture, restricted to the subspace of diffeomorphic
integrable $\omega$-compatible almost-K\"ahler structures,
gives many insights in the theory of extremal K\"ahler
metrics, where the leading conjectures are derived by a considerable
scope of analogy with the well-established GIT in the finite
dimensional case~\cite{fog-kir-mum}. It also suggests that
extremal almost-K\"ahler metrics would provide a natural extension of
the theory of extremal K\"ahler metrics to the non-integrable case. In
fact, this link has already become explicit in the toric case
\cite{don1}, where the existence of an extremal K\"ahler metric is
conjecturally equivalent to the existence of (infinitely many) non-integrable extremal almost-K\"ahler metrics;
this link was also used in \cite{apo-cal-gau-Fri 2} to find an {\it explicit} criterion to test K-stability (and therefore (non) existence of extremal K\"ahler metrics)
%this was further used in \cite{apo-cal-gau-Fri 2} to find an {\it explicit} criterion for the existence of  extremal K\"ahler metrics
on projective plane bundles over a curve and construct explicit examples of extremal almost-K\"ahler metrics. Besides the now appealing motivation, a systematic study of extremal almost-K\"ahler metric is still to come (see however \cite{apo-dra,kim-sun,leb-sim}).

In this paper, we generalize the notion of the Futaki invariant and extremal vector field to the general almost-K\"ahler case. This
amounts to the observation that fixing a compact subgroup $G\subset \textrm{Ham}(M,\omega)$ and considering $G$-invariant $\omega$-compatible almost-K\"ahler structures $(g,J)$, the $L^2$-projection of
the hermitian scalar curvature to the finite dimensional space of
hamiltonians of elements of $\textrm{Lie}(G)$ is independent of $(g,J)$. This
fact, which easily follows from the formal picture described above and
which is certainly known to experts (see e.g. \cite{apo-cal-gau-Fri 2,gau}), is established by a direct argument in Sec. \ref{sec3.2} below. By
taking $G$ be a maximal torus $T\subset \textrm{Ham}(M,\omega)$, we show in Lemma \ref{lem10} that the projection of the hermitian scalar curvature of any
$T$-invariant compatible almost-K\"ahler metric defines a Killing potential, $z_{\omega}^T$, which must coincide with the hermitian scalar curvature
of any $T$-invariant extremal almost-K\"ahler metric (should it exist). We
call the vector field $Z_{\omega}^T =grad_{\omega} z_{\omega}^T$ the
{\it extremal vector field} relative to $T.$ Its vanishing is an obstruction to the existence of $T$-invariant almost-K\"ahler metrics of constant hermitian scalar curvature. Noting that $Z_{\omega}^T$
doesn't change under a $T$-invariant isotopy of $\omega$, it
naturally generalizes the extremal vector field introduced
by \cite{fut-mab} in the K\"ahler context. 

The main technical input of the paper is proving that $Z_{\omega}^T$ has closed orbits when $[\frac{\omega}{2\pi}]$ is an integral class modulo torsion. This extends
the corresponding result in the K\"ahler case~\cite{nagak}, but we need
to adapt the `symplectic' approach of \cite{fut-mab1} which relies on the
localization formula~\cite{dui-heck}.
%, rather than the equivariant Riemann--Roch formula as in \cite{tian}.
Yet, the main technical ingredient is our Lemma \ref{lem4} which essentially computes the momentum map
of the $T$-action, with respect to the hermitian Ricci form $\rho^\nabla.$

The structure of the paper is as follows. In Sec. $2$, we introduce the necessary background of almost-K\"ahler geometry with special attention to holomorphic vector fields on almost-K\"ahler manifolds.
In particular, we obtain Bochner formulae involving the hermitian Ricci form $\rho^\nabla$ and the so-called $\ast$-Ricci form $\rho^\ast$ and derive some vanishing results (Corollaries \ref{corol3} and \ref{corol5}). 
We define in Sec. $3$ the extremal vector field relative to a maximal torus $T\subset \textrm{Ham}(M,\omega)$ in the general almost-K\"ahler case and prove its periodicity when $\left[\frac{\omega}{2\pi}\right]$ is integral modulo torsion. In Sec. $4$, we give an explicit (local) formula of the hermitian scalar curvature in Darboux coordinates which allows us to recast the expressions of the hermitian Ricci form and the scalar curvature in the toric case \cite{abr,don1}. We then specialise to the $4$-dimensional case and construct
infinite dimensional families of non-integrable extremal almost-K\"ahler metrics by using local toric symmetry; this allows us to obtain examples of non-integrable Hermitian-Einstein almost-K\"ahler metrics saturating LeBrun's estimates \cite{leb}.
Finally, in Sec. $5$, we discuss some avenues of further research.

\section{Preliminaries}
\subsection{Almost-K\"ahler structures}
An {\it almost-K\"ahler structure} on a real manifold $M$ of  dimension $2n$ is given by  a triple $(\omega,J,g)$ of a symplectic form $\omega$, an almost-complex structure $J$ and a Riemannian metric $g$, which satisfy the compatibility
relation
\begin{equation}\label{equa22}
g(\cdot,\cdot)=\omega(\cdot,J\cdot).
\end{equation}
We say that the almost-complex structure $J$ is $\omega${\it -compatible} if it induces a Riemannian metric via (\ref{equa22}).
If, additionally, the almost complex structure $J$ is {\it integrable}, then we have a {\it K\"ahler structure} on $M$.
%The complex vector bundle $(T(M),J)$ is identified with $T^{1,0}(M)$ via the map $X\to X^{1,0}=\frac{1}{2}\left(X-\sqrt{-1}JX\right)$ (so $X^{0,1}=\frac{1}{2}\left(X+\sqrt{-1}JX\right)$).
%Given an almost-complex structure $J$, the complexified tangent bundle splits as $T(M)\otimes\mathbb{C}=T^{1,0}(M)\oplus T^{0,1}(M)$ where $T^{1,0}(M)$, resp. $T^{0,1}(M)$ corresponds to the eigenvalue $\sqrt{-1}$, resp. $-\sqrt{-1}$
%under the $\mathbb{C}$-linear action of $J$.  The complex vector bundle $(T(M),J)$ is identified with $T^{1,0}(M)$ via the map $X\to X^{1,0}=\frac{1}{2}\left(X-\sqrt{-1}JX\right)$ (so $X^{0,1}=\frac{1}{2}\left(X+\sqrt{-1}JX\right)$).
%The almost-complex structure $J$ is integrable if and only if the distribution $T^{0,1}(M)$ is integrable. Also,
By the Newlander--Nirenberg theorem, the almost-complex structure $J$ is integrable if and only if the {\it Nijenhuis tensor} $N$
\begin{equation*}
4N(\cdot,\cdot)=\left[J\cdot,J\cdot\right]-\left[\cdot,\cdot\right]-J\left[J\cdot,\cdot\right]-J\left[\cdot,J\cdot\right]
\end{equation*}
vanishes; here $\left[\cdot,\cdot\right]$ stands for the Lie bracket. 

The complexified tangent bundle splits as 
\begin{equation*}
T(M)\otimes\mathbb{C}=T^{1,0}(M)\oplus T^{0,1}(M),
\end{equation*}
where $T^{1,0}(M)$ (resp. $T^{0,1}(M)$) corresponds to the eigenvalue $\sqrt{-1}$ (resp. $-\sqrt{-1}$) under the $\mathbb{C}$-linear action of $J$. The complex vector bundle $(T(M),J)$ is identified with $T^{1,0}(M)$ via the map $X\mapsto X^{1,0}=\frac{1}{2}\left(X-\sqrt{-1}JX\right)$ (resp. $X^{0,1}=\frac{1}{2}\left(X+\sqrt{-1}JX\right)$). The almost-complex structure $J$ also induces a decomposition of the complexified cotangent bundle 
\begin{equation*}
T^\ast(M)\otimes\mathbb{C}=\wedge^{1,0}(M)\oplus\wedge^{0,1}(M)
\end{equation*}
where $\wedge^{1,0}(M)$ (resp. $\wedge^{0,1}(M)$) is the annihilator of $T^{0,1}(M)$ (resp.  $T^{1,0}(M)$).  The almost-complex structure $J$ acts on the cotangent bundle $T^\ast(M)$ by $\left(J\alpha\right)(X)=-\alpha(JX)$. This action is extended to any section $\psi$ of the vector bundle $\wedge^p(M)$ of (real) $p$-forms by $J\psi(X_1,\cdots,X_p)=(-1)^p\psi(JX_1,\cdots,JX_p)$. 

Any section $\psi$ of $\otimes^2T^\ast(M)$ (and therefore of $T^\ast(M)\otimes T(M)$ which is identified to $\otimes^2T^\ast(M)$ via the metric) admits an orthogonal splitting $\psi=\psi^++\psi^-$, where $\psi^+$ is the $J$-invariant part and $\psi^-$ is the $J$-anti-invariant part,
given by
\begin{equation*}
\psi^+(\cdot,\cdot)=\frac{1}{2}\left(\psi(\cdot,\cdot)+\psi(J\cdot,J\cdot)\right){\text{ and }} \psi^-(\cdot,\cdot)=\frac{1}{2}\left(\psi(\cdot,\cdot)-\psi(J\cdot,J\cdot)\right).
\end{equation*}
In particular, the bundle of $2$-forms decomposes under the action of $J$
\begin{equation}\label{split1}
\wedge^2(M)=\mathbb{R}\,.\,\omega\oplus\wedge_0^{J,+}(M)\oplus\wedge^{J,-}(M),
\end{equation}
where $\wedge_0^{J,+}(M)$ is the subbundle of the {\it primitive} $J$-invariant $2$-forms (i.e. $2$-forms pointwise orthogonal to $\omega$) and $\wedge^{J,-}(M)$ is the subbundle of $J$-anti-invariant $2$-forms.

The fact that $\omega$ is closed implies the following identities (see \cite{kob-nom})
\begin{equation}\label{equa1}
g((D^g_XJ)Y,Z)+g((D^g_YJ)Z,X)+g((D^g_ZJ)X,Y)=0,
\end{equation}
\begin{equation}\label{equa2}
\left(D^g_X\omega\right)\left(Y,Z\right)=2g\left(JX,N(Y,Z)\right),
\end{equation}
where $D^g$ is the Levi-Civita connection with respect to the Riemannian metric $g$ and $X,Y,Z$ are any vector fields. Since $N$ is a $J$-anti-invariant $2$-form with values in $T(M)$, it follows from (\ref{equa2}) that
\begin{equation}\label{equa3}
D^g_{JX}J=-JD^g_XJ.
\end{equation}
Moreover, we readily deduce from the relation (\ref{equa2}) that the Nijenhuis tensor is identically zero if and only if $\omega$ (or equivalently $J$) is $D^g$-parallel.

From the exterior derivative $d$, we can define the twisted exterior differential $d^c=(-1)^pJdJ$ acting on $p$-forms (in particular $d^cf=Jdf$ for a smooth function $f$). A direct computation shows the following relation for any smooth function $f$ 
\begin{equation}\label{equa20}
dd^cf+d^cdf=4d^cf\left(N(\cdot,\cdot)\right).
\end{equation}
%More precisely,
%\begin{align*}
%dd^cf(X,Y)&=X.\left(d^cf(Y)\right)-Y.\left(d^cf(X)\right)-d^cf\left(\left[X,Y\right]\right)\\
%&=-X.(JY.f)+Y(JX.f)+J\left[X,Y\right].f,
%\end{align*}
%whereas
%\begin{align*}
%d^cdf(X,Y)&=-dJdf(JX,JY)\\
%&=-JX.(Y.f)+JY(X.f)-J\left[JX,JY\right].f.
%\end{align*}
It follows that the almost-complex structure $J$ is integrable if and only if $d$ and $d^c$ anticommute.

We denote by $\Lambda_\omega$ the contraction by the symplectic form $\omega$, defined for a $p$-form $\psi$ by $\Lambda_{\omega}(\psi)=\frac{1}{2}\sum_{i=1}^{2n}\psi\left(e_i,Je_i,\cdots\right),$ where $\{e_i\}$ is a local $J$-adapted orhonormal frame. As noticed by Gauduchon \cite{gau} and Merkulov \cite{mer}, the commutator of $\Lambda_\omega$ and $d$ is equal to
\begin{equation}\label{equa15}
\left[\Lambda_\omega,d\right]=-\delta^c,
\end{equation}
where $\delta^c=(-1)^pJ\delta^g J$ is the twisted codifferential acting on $p$-forms; here $\delta^g$ is the codifferential defined as the formal adjoint of $d$ with respect to the metric $g$ ($\delta^g$ also stands for the adjoint of $D^g$ with respect to $g$ when it is applied to sections of $\otimes^pT^\ast(M)$).
The operator $\Lambda_\omega$ commutes with $J$, so the relation (\ref{equa15}) implies
\begin{equation}\label{equa16}
\left[\Lambda_\omega,d^c\right]=\delta^g.
\end{equation}
It follows from (\ref{equa16}) that on any almost-K\"{a}hler manifold we have \cite{gau}
\begin{equation}\label{equa17}
\delta^g d^c+d^c\delta^g=0.
\end{equation}
\subsection{The hermitian connection}
The canonical hermitian connection $\nabla$ corresponding to $J$ is defined by
\begin{equation*}
\nabla_XY=D^g_XY-\frac{1}{2}J\left(D^g_XJ\right)Y.
\end{equation*}
The connection $\nabla$ preserves $J$ (i.e. $\nabla J=0$) and has a $J$-anti-invariant torsion.
It also induces the Cauchy--Riemann operator  $\bar{\partial}$ on $T^{1,0}(M)$ \cite{gau}, where we recall $({\bar{\partial}}Y^{1,0})({X^{0,1}} ):= \left[X^{0,1},Y^{1,0}\right]^{1,0}.$
Indeed,
\begin{proposition}\label{prop5}
For any vector fields $X,Y$
\begin{equation*}
\nabla_{X^{0,1}}Y^{1,0}=\left[X^{0,1},Y^{1,0}\right]^{1,0}. 
\end{equation*}
\end{proposition}
\begin{proof}
Using the above defintion of $\nabla$ and the relation (\ref{equa3}), we have
\begin{align*}
\nabla_{X^{0,1}}Y^{1,0}&=D^g_{X^{0,1}}Y^{1,0}-\frac{1}{2}J\left(D^g_{X^{0,1}}J\right)Y^{1,0}\\
&=\frac{1}{4}D^g_{X+\sqrt{-1}JX}\left(Y-\sqrt{-1}JY\right)-\frac{1}{8}J\left(D^g_{X+\sqrt{-1}JX}J\right)\left(Y-\sqrt{-1}JY\right)\\
&=\frac{1}{4}\Big{(}D^g_XY+D^g_{JX}(JY)\Big{)}+\frac{1}{4}\sqrt{-1}\Big{(}D^g_{JX}Y-D^g_X{(JY)}\Big{)}\\
&=\frac{1}{4}\left(Id-\sqrt{-1}J\right)\Big{(}D^g_XY+D^g_{JX}(JY)\Big{)}.
\end{align*}
On the other hand,
\begin{align*}
\left[X^{0,1},Y^{1,0}\right]^{1,0}&=\frac{1}{4}\left[X+\sqrt{-1}JX,Y-\sqrt{-1}JY\right]^{1,0}\\
&=\frac{1}{8}\left(Id-\sqrt{-1}J\right)\left[X+\sqrt{-1}JX,Y-\sqrt{-1}JY\right]\\
&=\frac{1}{8}\left(Id-\sqrt{-1}J\right)\Big{(}\left[X,Y\right]+\left[JX,JY\right]+J\left[JX,Y\right]-J\left[X,JY\right]\Big{)}\\
&=\frac{1}{8}\left(Id-\sqrt{-1}J\right)\Big{(}D^g_XY-D^g_YX+D^g_{JX}(JY)-D^g_{JY}(JX)\\
&+JD^g_{JX}Y-JD^g_Y(JX)-JD^g_X(JY)+JD^g_{JY}X\Big{)}\\
&=\frac{1}{8}\left(Id-\sqrt{-1}J\right)\Big{(}D^g_XY-D^g_YX+D^g_{JX}(JY)+J(D^g_{Y}J)X-J(D^g_{JY}X)\\
&+JD^g_{JX}Y+D^g_YX-J(D^g_YJ)X+D^g_XY+(D^g_{JX}J)Y+JD^g_{JY}X\Big{)}\\
&=\frac{1}{4}\left(Id-\sqrt{-1}J\right)\Big{(}D^g_XY+D^g_{JX}(JY)\Big{)}.
\end{align*}
The proposition follows.
\end{proof}
\subsection{Ricci forms}\label{ricci}
The canonical hermitian connection $\nabla$ on $T(M)$ induces a hermitian connection on the anti-canonical bundle $K_J^{-1}(M)=\wedge^{n,0}(M)$ (equipped with the hermitian structure induced by $g$) with curvature $\sqrt{-1}\rho^{\nabla}$, where $\rho^{\nabla}$ is a closed (real) $2$-form, called the {\it hermitian Ricci form}. Hence, the $2$-form $\rho^{\nabla}$ is a deRham representative of $2\pi c_1(M,J)$ in $H^2(M,\mathbb{R})$ where $ c_1(M,J)$ is the (real) Chern class of $K_J^{-1}(M)$. 

We consider also $\rho^\ast=R(\omega),$ the image of the symplectic form $\omega$ by the (Riemannian) curvature operator $R$. The $2$-form $\rho^\ast$ is called the $\ast${\it -Ricci form}. %In the same way, we define the $\ast${\it -scalar curvature} $s^\ast$
%given by
%\begin{equation*}
%s^\ast\omega^n=2n\left(\rho^\ast\wedge\omega^{n-1}\right).
%\end{equation*}
We have the following relation between these Ricci-type tensors (see \cite{apo-dra})
\begin{equation}\label{equa12}
\rho^{\nabla}(X,Y)=\rho^\ast(X,Y)-\frac{1}{4}tr(JD^g_XJ\circ D^g_YJ).
\end{equation}
In the K\"ahler case (i.e. when $D^g\omega=D^gJ=0$), we readily deduce from the relation (\ref{equa12}) that $\rho^{\nabla}=\rho^\ast,$ which is also equal to the usual Ricci form. Note that neither $\rho^\nabla$ nor
$\rho^\ast$ is $J$-invariant in general. In fact, by (\ref{equa12}) and (\ref{equa3}), $\rho^\nabla$ is $J$-invariant if and only if $\rho^\ast$ is.

%An almost- K\"ahler metric is called {\it hermitian-Einstein} if the hermitian Ricci form $\rho^{\nabla}$ is a (constant) multiple of the symplectic form, i.e. $\rho^{\nabla}=\frac{s^{\nabla^g}}{2n}\omega.$
\subsection{Holomorphic vectors fields on an almost-K{\"a}hler manifold}
In this section, we study (pseudo-)holomorphic vector fields on an almost-K{\"a}hler manifold $(M^{2n},\omega,J,g)$.
A (real) vector field $X$ is said {\it holomorphic} if $\mathfrak{L}_{X}J=0,$ where $\mathfrak{L}$ denotes the Lie derivative.
It is equivalent to say that $\left[X,JY\right]=J\left[X,Y\right]$ for any vector field $Y.$
\begin{lemma}\label{lem2}
On an almost-K{\"a}hler manifold $(M^{2n},\omega,J,g)$, $X$ is a holomorphic vector field if and only if
\begin{equation}\label{equa6}
{\left(D^{g}X\right)}^-=-\frac{1}{2}D^{g}_{JX}J.
\end{equation}
\end{lemma}
\begin{proof}
For any vector field $X$, we have $\mathfrak{L}_{X}J=D^{g}_{X}J-\left[D^{g}X,J\right]$. In
particular, if $X$ is holomorphic then $D^{g}_{X}J=\left[D^{g}X,J\right]$. On the other hand
$J\left[D^{g}X,J\right]=2{\left(D^{g}X\right)}^-$. The equality (\ref{equa6}) followns from  (\ref{equa3}).
\end{proof}
\begin{lemma}\label{lem3}
Let $\alpha$ be a $1$-form. We have
\begin{equation*}
\delta^g {\left(D^g\alpha\right)}^+-\delta^g{\left(D^g\alpha\right)}^-=\rho^\ast\left({\alpha^{\sharp},J\cdot}\right)-\sum\limits_{i=1}^{2n}(D^g_{Je_i}\alpha)\left((D^g_{e_i}J)(\cdot)\right),
\end{equation*}
where $\{e_i\}$ is a local $J$-adapted orhonormal frame, $\delta^g$ the formal adjoint of the Levi-Civita connection $D^g$ with respect to the metric $g$ and $\sharp$ stands for the isomorphism between $T^\ast(M)$ and $T(M)$ induced by $g^{-1}.$
\end{lemma}
\begin{proof}
By using the fact that $\delta^g J=-\sum\limits_{i=1}^{2n}\left(D^g_{e_i}J\right)(e_i)=0$, we have
\begin{align*}
\left(\delta^g {\left(D^g\alpha\right)}^+-\delta^g{\left(D^g\alpha\right)}^-\right)(X)&=-\sum\limits_{i=1}^{2n}D^g_{e_i}\left({\left(D^g\alpha\right)}^+-{\left(D^g\alpha\right)}^-\right)(e_i,X)\\
&=-\sum\limits_{i=1}^{2n} \left[ D^g_{e_i} \left[(D^g_{Je_i}\alpha)(JX) \right]-(D^g_{D^g_{e_i}(Je_i)}\alpha)(JX)-(D^g_{Je_i}\alpha)(JD^g_{e_i}X)\right]
\end{align*}
\begin{align*}
\;\;\;\;\;\:\;\;\;\;\;\;\:\;\;\;\;\;\;\:\;\;\;\;\;\;\;\;\;\:\;\;\;\;\;\;\;\;\;\:\;\;\;\;\;\;\;\,&=-\sum\limits_{i=1}^{2n}\bigg{[}\Big{(}D^g_{e_i}(D^g_{Je_i}\alpha)\Big{)}(JX)+(D^g_{Je_i}\alpha)\Big{(}D^g_{e_i}(JX)\Big{)}-(D^g_{D^g_{e_i}(Je_i)}\alpha)(JX)\\
&-(D^g_{Je_i}\alpha)(JD^g_{e_i}X)\bigg{]}\\
&=-\sum\limits_{i=1}^{2n} \bigg{[}\Big{(}D^g_{e_i}(D^g_{Je_i}\alpha)\Big{)}(JX)-(D^g_{D^g_{e_i}(Je_i)}\alpha)(JX)+(D^g_{Je_i}\alpha)\Big{(}(D^g_{e_i}J)(X)\Big{)}\bigg{]}  \\
&=-\sum\limits_{i=1}^{2n}(D^2_{e_i,Je_i}\alpha)(JX)-\sum\limits_{i=1}^{2n}(D^g_{Je_i}\alpha)\Big{(}(D^g_{e_i}J)(X)\Big{)}\\
&=\frac{1}{2}\sum\limits_{i=1}^{2n}(R_{e_i,Je_i}\alpha)(JX)-\sum\limits_{i=1}^{2n}(D^g_{Je_i}\alpha)\Big{(}(D^g_{e_i}J)(X)\Big{)}\\
&=\rho^\ast_{\alpha^{\sharp},JX}-\sum\limits_{i=1}^{2n}(D^g_{Je_i}\alpha)\Big{(}(D^g_{e_i}J)(X)\Big{)}.
\end{align*}
\end{proof}
\begin{corollary}\label{corol2}
Let $\alpha$ be a $1$-form such that $X=\alpha^\sharp$ is holomorphic vector field. Then
\begin{equation*}
\delta^g {\left(D^g\alpha\right)}^+-\delta^g{\left(D^g\alpha\right)}^-=\rho^{\nabla}\left({X,J\cdot}\right).
\end{equation*}
\end{corollary}
\begin{proof}
Combining Lemmas \ref{lem2} and \ref{lem3} and using the relations (\ref{equa1}) and (\ref{equa12}), we have
\begin{align*}
\left(\delta^g {\left(D^g\alpha\right)}^+-\delta^g{\left(D^g\alpha\right)}^-\right)(Y)&=\rho^\ast_{X,JY}-\sum\limits_{i=1}^{n}(D^g_{Je_i}\alpha)\Big{(}(D^g_{e_i}J)(Y)\Big{)}\\
&=\rho^{\nabla}_{X,JY}+\frac{1}{4}tr(JD^g_XJ\circ D^g_{JY}J)-\sum\limits_{i=1}^{2n}(D^g_{Je_i}\alpha)\Big{(}(D^g_{e_i}J)(Y)\Big{)}
\end{align*}
\begin{align*}
\;\;\;\;\;\:\;\;\;\;\;\;\:\;\;\;\;\;\;\:\;\;\;\;\;\;\;\;\;\:\;\;\;\;\;\;\;\;\;\:\;\;\;\;\;\;\,\,&=\rho^{\nabla}_{X,JY}+\frac{1}{4}tr(JD^g_XJ\circ D^g_{JY}J)-\sum\limits_{i=1}^{2n}\left({(D^g\alpha)}^-_{Je_i}\right)\Big{(}(D^g_{e_i}J)(Y)\Big{)}\\
&=\rho^{\nabla}_{X,JY}+\frac{1}{4}tr(JD^g_XJ\circ D^g_{JY}J)+\frac{1}{2}\sum\limits_{i=1}^{2n}g\Big{(}(D^g_{JX}J)(Je_i),(D^g_{e_i}J)(Y)\Big{)}\\
&=\rho^{\nabla}_{X,JY}+\frac{1}{4}tr(JD^g_XJ\circ D^g_{JY}J)-\frac{1}{2}\sum\limits_{i=1}^{2n}g\Big{(}(D^g_{X}J)(e_i),(D^g_{e_i}J)(Y)\Big{)}\\
&=\rho^{\nabla}_{X,JY}+\frac{1}{4}tr(JD^g_XJ\circ D^g_{JY}J)-\frac{1}{2}\sum\limits_{i,k=1}^{2n}\Big{(}(D^g_{X}J)(e_i,e_k)\otimes(D^g_{e_i}J)(Y,e_k)\Big{)}\\
&=\rho^{\nabla}_{X,JY}+\frac{1}{4}tr(JD^g_XJ\circ D^g_{JY}J)\\
&-\frac{1}{4}\sum\limits_{i,k=1}^{2n}\bigg{(}(D^g_{X}J)(e_i,e_k)\otimes\Big{[}(D^g_{e_i}J)(Y,e_k)-(D^g_{e_k}J)(Y,e_i)\Big{]}\bigg{)}\\
&=\rho^{\nabla}_{X,JY}+\frac{1}{4}tr(JD^g_XJ\circ D^g_{JY}J)-\frac{1}{4}\sum\limits_{i,k=1}^{2n}\Big{(}(D^g_{X}J)(e_i,e_k)\otimes(D^g_{Y}J)(e_i,e_k)\Big{)}\\
&=\rho^{\nabla}_{X,JY}+\frac{1}{4}tr(JD^g_XJ\circ D^g_{JY}J)-\frac{1}{4}\sum\limits_{i=1}^{2n}g\Big{(}(D^g_{X}J)(e_i),(D^g_{Y}J)(e_i)\Big{)}\\
&=\rho^{\nabla}_{X,JY}-\frac{1}{4}tr(D^g_XJ\circ D^g_{Y}J)+\frac{1}{4}\sum\limits_{i=1}^{2n}g\Big{(}(D^g_{X}J)(D^g_{Y}J)e_i,e_i\Big{)}\\
&=\rho^{\nabla}_{X,JY}-\frac{1}{4}tr(D^g_XJ\circ D^g_{Y}J)+\frac{1}{4}tr(D^g_XJ\circ D^g_{Y}J)\\
&=\rho^{\nabla}_{X,JY}.
\end{align*}
\end{proof}
\begin{corollary}\label{corol3}
Let $(M^{2n},\omega,J,g)$ be a compact almost-K\"ahler manifold. Suppose that the tensor $(\rho^\ast)^+(\cdot,J\cdot)$ is negative-definite. Then, there is no non-trivial holomorphic vector field on $M$.
\end{corollary}
\begin{proof}
Let $\alpha$ be the dual (by the metric $g$) of a holomorphic vector field $X$. By Lemma~\ref{lem3}, we have
\begin{equation}\label{equa9}
\left(\delta^g {\left(D^g\alpha\right)}^+-\delta^g{\left(D^g\alpha\right)}^-\right)(X)=\rho^\ast(X,JX)-\sum\limits_{i=1}^{2n}(D^g_{Je_i}\alpha)\left((D^g_{e_i}J)(X)\right).
\end{equation}
Using Lemma \ref{lem2}, we simplify the second term of the right hand side of (\ref{equa9}) as in the proof of Corollary \ref{corol2}
\begin{align*}
\sum\limits_{i=1}^{n}g\Big{(}(D^g_{Je_i}\alpha),(D^g_{e_i}J)(X)\Big{)}&=\sum\limits_{i=1}^{2n}g\Big{(}(D^g\alpha)^-_{Je_i},(D^g_{e_i}J)(X)\Big{)}\\
&=-\frac{1}{2}\sum\limits_{i=1}^{2n}g\Big{(}(D^g_{JX}J)(Je_i),(D^g_{e_i}J)(X)\Big{)}
\end{align*}
\begin{align*}
\;\;\;\;\;\:\;\;\;\;\;\;\:\;\;\;\;\;\;\:\;\;\;\;\;\;\;\;\;\:\;\;\;\;\;\;\;\;\;\:\;\;\;\;\;\;\,\,
&=\frac{1}{2}\sum\limits_{i=1}^{2n}g\Big{(}(D^g_{X}J)(e_i),(D^g_{e_i}J)(X)\Big{)}\\
&=\frac{1}{2}\sum\limits_{i,k=1}^{2n}\Big{(}(D^g_{X}J)(e_i,e_k)\otimes(D^g_{e_i}J)(X,e_k)\Big{)}\\
&=\frac{1}{4}\sum\limits_{i,k=1}^{2n}\Big{(}(D^g_{X}J)(e_i,e_k)\otimes\left[(D^g_{e_i}J)(X,e_k)-(D^g_{e_k}J)(X,e_i)\right]\Big{)}\\
&=\frac{1}{4}\sum\limits_{i,k=1}^{2n}\Big{(}(D^g_{X}J)(e_i,e_k)\otimes(D^g_{X}J)(e_i,e_k)\Big{)}\\
&=\frac{1}{4}\sum\limits_{i=1}^{2n}\Big{(}(D^g_{X}J)(e_i),(D^g_{X}J)(e_i)\Big{)}.
\end{align*}
Integrating (\ref{equa9}), we obtain
\begin{equation*}
\left|\left(D^g\alpha\right)^+\right|_{L^2}^2-\left|\left(D^g\alpha\right)^-\right|_{L^2}^2=\left(\int_M\rho^\ast(X,JX)\frac{\omega^n}{n!}\right)-\frac{1}{4}\left|D^g_{X}J\right|_{L^2}^2.
\end{equation*}
Because of Lemma \ref{lem2}, $\left|\left(D^g\alpha\right)^-\right|_{L^2}^2=\frac{1}{4}\left|D^g_{X}J\right|_{L^2}^2$. Then, $\left|\left(D^g\alpha\right)^+\right|_{L^2}^2=\Big{(}\int_M\rho^\ast(X,JX)\frac{\omega^n}{n!}\Big{)}$. The corollary follows
\end{proof}
%\begin{corollary}\label{corol4}
%Let $(M^{2n},\omega,J,g)$ be a compact almost-K\"ahler
%manifold. Suppose that $(\rho^{\nabla})^+(X,JX)<0$ for any vector field $X$. Then, there is no holomorphic vector field $X$ such that $JX$ is holomorphic too.
%\end{corollary}
%\begin{proof}
%A direct calculation shows that $\mathfrak{L}_{JX}J-J\mathfrak{L}_{X}J=4N(X,\cdot)$ for any vector field $X$. In particular, If $X$ and $JX$ are holomorphic, we deduce that $N(X,\cdot)=0$ and then, by the relations (\ref{equa1}) and (\ref{equa2}), we  have
%$D^g_XJ=0$. Using the formula (\ref{equa12}), we obtain
%$\rho^{\nabla}(X,JX)=\rho^\ast (X,JX)$. Let $\alpha$ be the dual of $X$ by the metric $g.$ By the same arguments in the proof of corollary~\ref{corol3}, we have
%$\left|\left(D^g\alpha\right)^+\right|^2=\left(\int_M\rho^\ast(X,JX)\frac{\omega^n}{n!}\right)=\left(\int_M\rho^{\nabla}(X,JX)\frac{\omega^n}{n!}\right)$. The corollary follows
%\end{proof}
\begin{lemma}\label{lem4}
Let $(M^{2n},\omega,J,g)$ be an almost-K\"ahler
manifold. If $X$ is a hamiltonian Killing vector field, then
\begin{equation}\label{equa11}
-\frac{1}{2}d\Delta^g f^X =\rho^{\nabla}(X,\cdot),
\end{equation}
where $f^X$ is the momentum of $X$ with respect to $\omega$ (i.e. $\omega(X,\cdot)=-df^X$) and $\Delta^g$ is the Riemannian Laplacian with respect to $g$.
\end{lemma}
\begin{proof}
By hypothesis, $X=\left(d^cf\right)^{\sharp}$ is a Killing vector field then $D^gd^cf=\frac{1}{2}dd^cf$. So, we have $\delta^g D^gd^cf=\delta^g\left(D^gd^cf\right)^++\delta^g\left(D^gd^cf\right)^-=\frac{1}{2}\delta^g dd^cf$. Combining this with Corollary~\ref{corol2} applied for the holomorphic vector field $X$, we obtain 
\begin{equation}\label{equa19}
2\delta^g \left(D^gd^cf\right)^-=\frac{1}{2}\delta^g dd^cf-\rho^{\nabla}({X,J\cdot}).
\end{equation}
On the other hand (by Lemma \ref{lem2})
$\left(D^gd^cf\right)^-=\frac{1}{2}D^g_{(df)^\sharp}\omega$. Then, using the relations (\ref{equa2}) and (\ref{equa20}), we obtain $\left(D^gd^cf\right)^-= \frac{1}{4}(dd^cf+d^cdf)$ and therefore
\begin{equation}\label{equa18}
\delta^g\left(D^gd^cf\right)^-=\frac{1}{4}\delta^g(dd^cf+d^cdf).
\end{equation}
Combining (\ref{equa18}) with (\ref{equa19}), we obtain, using the relation (\ref{equa17})
\begin{eqnarray*}
\frac{1}{2}\delta^g(dd^cf+d^cdf)&=&\frac{1}{2}\delta^g dd^cf-\rho^{\nabla}({X,J\cdot})\\
\frac{1}{2}\delta^g d^cdf &=&-\rho^{\nabla}({X,J\cdot})\\
-\frac{1}{2}d^c\delta^g df &=&-\rho^{\nabla}({X,J\cdot})\\
-\frac{1}{2}d\Delta^g f &=&\rho^{\nabla}({X,\cdot}).
\end{eqnarray*}
\end{proof}
\begin{corollary}\label{corol5}
Let $(M^{2n},\omega,J,g)$ be a compact almost-K\"ahler
manifold. Suppose that the tensor $(\rho^{\nabla})^{+}(\cdot,J\cdot)$ is negative-semidefinite. Then, there is no non-trivial hamiltonian Killing vector field on $M$.
\end{corollary}
\begin{proof}
Suppose that $X=\left(d^cf\right)^{\sharp}$ is a hamiltonian Killing vector field. By Lemma~\ref{lem4}, we have
\begin{eqnarray*}
-\frac{1}{2}d\Delta^g f\left((df)^\sharp\right) =\rho^{\nabla}\left(\left(d^cf\right)^{\sharp},(df)^\sharp\right).
\end{eqnarray*}
By integrating, we obtain $\frac{1}{2}\left|\Delta^gf\right|_{L^2}^2=\left(\int_M\rho^{\nabla}\Big{(}\left(df\right)^{\sharp},(d^cf)^\sharp\right)\frac{\omega^n}{n!}\Big{)}$.
If $\Delta^gf$ is identically zero on a compact Riemannian manifold, $f$ must be constant and thus $X=0$. The corollary follows.
\end{proof}

\subsection{A localization formula}
Let $(M^{2n},\omega)$ be a compact symplectic manifold endowed with a hamiltonian $S^1$-action. Let $X$ be the generator of this action with a momentum $f^X$ i.e. $\omega(X,\cdot)=-df^X.$
Thus, $X$ is a hamiltonian Killing vector field with respect to some compatible almost-K\"{a}hler metrics.
The fixed points of the action form a finite union of connected symplectic submanifolds \cite{duf-sal} $N_1,\cdots,N_\gamma$ such that $f^X$ is constant on each $N_j$. For each $N_j$, the normal bundle $E_j$ splits into complex line bundles
$E_j=L_1^j\oplus\cdots\oplus L_{m_j}^j$ on which $S^1$ acts with integer weights $k_1^j,\cdots,k_{m_j}^j$. Then, we have the following  formula (see \cite{dui-heck,duf-sal}) 
\begin{equation}\label{equa13}
\int_Me^{-hf^X}\frac{{\omega}^n}{n!}=\sum_{j=0}^{\gamma}\left(\int_{N_j}\prod_{i=1}^{m_j}\frac{1}{c_1(L_i^j)+k_i^jh}\right)e^{-hf^X(N_j)},
\end{equation}
for every $h\in\mathbb{C}.$ Here, $c_1(\cdot)$ denote the first Chern class and we take the formal inverse
\begin{equation*}
\frac{1}{c_1(L_i^j)+k_i^jh}=\frac{1}{k_i^jh}\sum_{r=1}^{(dim N_j)/2}\left(-\frac{c_1(L_i^j)}{k_i^jh}\right)^r.
\end{equation*}

\section{Extremal almost-K\"ahler metrics}
For the rest of the paper, $(M^{2n},\omega)$ is a {\it compact} and {\it connected} symplectic manifold. We denote by $\left[\omega\right]$ the deRham cohomology class of $\omega.$ Any $\omega$-compatible almost-complex structure  is identified with the induced Riemannian metric.
\subsection{Hermitian scalar curvature as a moment map.}\label{sec3.1}

We define the {\it hermitian scalar curvature} $s^{\nabla}$ of an almost-K\"{a}hler structure $(\omega,J,g)$ as  the trace of $\rho^\nabla$ with respect to $\omega$, i.e.
\begin{equation}\label{equa21}
s^{\nabla}\omega^n=2n\left(\rho^{\nabla}\wedge\omega^{n-1}\right),
\end{equation}
or, in equivalent way,
\begin{equation}\label{equa24}
s^{\nabla}=2\Lambda_\omega\rho^{\nabla}.
\end{equation}

Denote by $AK_{\omega}$ the Fr\'echet space of $\omega$-compatible almost-complex structures.
The space $AK_{\omega}$ comes naturally equipped with a formal K\"{a}hler structure $(\mathbf{\Omega},\mathbf{J})$ (described first by Fujiki in \cite {fuj}). More precisely, the tangent space $T_J(AK_{\omega})$ at a point $J$ is identified with the space of
$g$-symmetric, $J$-anti-invariant fields of endomorphisms of $T(M)$ (where $(\omega,J,g)$ is the corresponding almost-K\"ahler metric). Then, for $A,B\in T_J(AK_{\omega})$, the K\"ahler form $\mathbf{\Omega}$ is given by $\mathbf{\Omega}_J(A,B)=\int_Mtr\left(J\circ A\circ B\right)\frac{\omega^n}{n!}$ while the $\mathbf{\Omega}$-compatible (integrable) almost-complex structure $\mathbf{J}$ is defined by $\mathbf{J}_JX=J\circ X.$

Let $\textrm{Ham}(M,\omega)$ be the group of hamiltonian symplectomorphisms of $(M^{2n},\omega)$. The Lie algebra of $\textrm{Ham}(M,\omega)$ is identified with the space of smooth functions on $M$ with zero mean value and the {\it Poisson bracket};
it is also equipped with an equivariant inner product, given by the $L^2$-norm with respect to $\frac{\omega^n}{n!}$

A key observation, made by Fujiki \cite {fuj} in the integrable case and by Donaldson \cite{don} in the general almost-K\"{a}hler case, asserts that the natural action of $\textrm{Ham}(M,\omega)$ on $AK_{\omega}$ is hamiltonian with momentum given by the hermitian scalar curvature $s^\nabla$. More precisely, the moment map is
\begin{equation}\label{equa35}
\mu_J(f)=\int_M s^{\nabla}f\frac{\omega^n}{n!}
\end{equation}
where $s^{\nabla}$ is the hermitian scalar curvature of the induced almost-K\"{a}hler metric $(\omega,J,g)$. The square-norm of the hermitian
scalar curvature defines a functional on $AK_{\omega}$ 
\begin{equation}\label{equa4}
J\mapsto\int_M {\left(s^{\nabla}\right)}^2\frac{\omega^n}{n!}.
\end{equation}
\begin{definition}
The critical points of the functional (\ref{equa4}) are called {\it extremal almost-K\"{a}hler metrics}. 
\end{definition}
%Hence, by the moment map approach, Donaldson extends naturally the Calabi problem  \cite{cal} to the non-integrable case.
The functional (\ref{equa4}) corresponds to the {\it square-norm function} of the moment map. This observation was used  \cite{apo-dra} to deduce that the critical points of (\ref{equa4}) are precisely the almost-K\"{a}hler metrics $g$ for which the symplectic gradient of their hermitian scalar curvature $grad_{\omega}s^{\nabla}$, defined as the symplectic dual of $ds^\nabla$ (i.e. $-ds^\nabla=\omega(grad_{\omega}s^{\nabla},\cdot)=g(Jgrad_{\omega}s^{\nabla},\cdot)$ and thus $grad_{\omega}s^{\nabla}=Jgrad\,s^\nabla$), is a Killing vector field of $g$; as it is hamiltonian, this is also equivalent to being holomorphic with respect to $J.$ Indeed, it follows from \cite{fog-kir-mum} that a point $x_0$ is critical for the square norm
function of the momentum map if and only if the image of $x_0$ by the moment map belongs to the Lie subalgebra corresponding to the stabilizer of $x_0$ by the action (where the Lie algebra is identified with its dual vector space via the inner product). In our context, this precisely means that $grad_{\omega}s^{\nabla}$ is a Killing hamiltonian vector field.

\begin{proposition}
A metric $g$ is a critical point of (\ref{equa4}) if and only if $grad_{\omega}s^{\nabla}$ is a Killing vector field with respect to $g.$
\end{proposition}
\begin{proof}
We reproduce here a direct verification made by Gauduchon \cite{gau} using Mohsen formula \cite{moh}.
The Mohsen formula states that, for a path $J_t\in AK_{\omega}$, the first variation of the connection $1$-forms $\alpha_t$ of the hermtian connections on $K_{J_t}^{-1}(M)$, induced by the canonical hermitian connections $\nabla_t$ corresponding to $J_t$, is given by
\begin{equation*}
\frac{d}{dt}\alpha_t=\frac{1}{2}\delta^{g_t}\dot{J},
\end{equation*}
where $\delta^{g_t}$ is the codifferential with respect to the metric $g_t(\cdot,\cdot)=\omega(\cdot,J_t\cdot)$ and $\dot{J}=\frac{d}{dt}J_t.$ Therefore, by definition, $\frac{d}{dt}\rho^{\nabla_t}=-\frac{1}{2}d\delta^{g_t}\dot{J}$. Hence, by (\ref{equa24}) and (\ref{equa15}), we obtain
$\frac{d}{dt}s^{\nabla_t}=-\Lambda_{\omega} d\delta^{g_t}\dot{J}=\delta_t^c\delta^{g_t}\dot{J}=-\delta^{g_t} J_t\delta^{g_t}\dot{J};$ here $\delta_t^c$ is the twisted codifferential with respect to $J_t.$
Therefore,
\begin{equation}\label{equa5}
\frac{d}{dt}s^{\nabla_t}=-\delta^{g_t} J_t\delta^{g_t}\dot{J}.
\end{equation}
Using (\ref{equa5}), we compute the first derivative of (\ref{equa4}) in the direction of $\dot{J}$
\begin{equation*}
\frac{d}{dt}\left( \int_M {\left(s^{\nabla_t}\right)}^2\frac{\omega^n}{n!}
\right)=2\int_Ms^{\nabla_t}\left(-\delta^{g_t} J_t\delta^{g_t}
\dot{J}\right)\frac{\omega^n}{n!}=2\int_Mg_t\left(D^{g_t}d^c_ts^{\nabla_t},\dot{J}\right)\frac{\omega^n}{n!},
\end{equation*}
where $d_t^c$ is the twisted exterior differential corresponding to $J_t$.
%From the identities $J_t^2=-I$ and $\omega(J_t\cdot,J_t\cdot)=\omega(\cdot,\cdot)$,  we deduce
However, $\dot{J}=\frac{d}{dt}J_t$ is a $g_t$-symmetric, $J_t$-anti-invariant endomorphism of $T(M)$. Hence, $J$ is a critical point if and only if
the symmetric, $J$-anti-invariant part of $D^{g}d^cs^{\nabla}$ is identically zero. On the other hand, for any vector field $X$ preserving $\omega$, we have
%Gauduchon shows by a direct calculation that the symmetric part of $D^gX$, for any vector field $X$ preserving the symplectic form $\omega$, is given by $-\frac{1}{2}J\mathfrak{L}_XJ$ and therefore $J$-anti-invariant. Indeed,
\begin{eqnarray*}
0&=& \left(\mathfrak{L}_X\omega\right)(JY,Z)\\
&=& -\left(\mathfrak{L}_Xg\right)(Y,Z)+g\left( \left(\mathfrak{L}_XJ\right)JY,Z\right)\\
&=&-g\left(\left(D^gX\right)Y,Z\right)-g\left(Y,\left(D^gX\right)Z\right)+g\left( \left(\mathfrak{L}_XJ\right)JY,Z\right).
\end{eqnarray*}
It follows that the symmetric part of $D^{g}d^cs^{\nabla}$ is already $J$-anti-invariant. Hence, $J$ is a critical point if and only if $D^{g}d^cs^{\nabla}$
is skew-symmetric which means that the symplectic gradient $grad_{\omega}s^{\nabla}$ is a Killing vector field.
\end{proof}

\begin{remark}
On a complex manifold $(M,J)$, the Calabi problem \cite{cal} consist in studying the Calabi functional given by the $L^2$-norm of the scalar curvature of the K\"ahler metrics whose
K\"ahler form belongs to a fixed K\"ahler class $\Omega=\left[\omega\right]$. 
It turns out that the critical points of the Calabi functional, called {\it Calabi extremal K\"ahler metrics}, are the K\"ahler metrics for which the symplectic gradient of the scalar curvature is a Killing vector field.
The extremal almost-K\"{a}hler metrics thus appear as a natural extension of the Calabi extremal K\"ahler metrics to the non-integrable case. Indeed, since any two K\"ahler forms in a fixed K\"ahler class $\Omega$ are isotopic, the K\"ahler metrics in $\Omega=\left[\omega\right]$ are embedded, via Moser's lemma \cite{duf-sal}, in the space of $\omega$-compatible integrable almost-complex structures $K_{\omega}.$ 
\end{remark}

\subsection{The extremal vector field}\label{sec3.2}

We fix a compact group $G$ in the (infinite dimensional) group $\textrm{Ham}(M,\omega)$ of hamiltonian symplectomorphisms of $(M^{2n},\omega)$.
Denote by $\mathfrak{g}_{\omega}\subset C^\infty(M)$ the finite dimensional
space of smooth functions which are hamiltonians with zero mean value of elements of
$\mathfrak{g}=Lie(G).$ It is well-known that $\mathfrak{g}_{\omega}$ has a Lie algebra structure given by the Poisson bracket and $\mathfrak{g}_{\omega}$ is isomorphic to $\mathfrak{g}.$
Denote by $\Pi_{\omega}$ the $L^2$-orthogonal
projection of a smooth function on $\mathfrak{g}_{\omega}$ with
respect to the volume form $\frac{\omega^n}{n!}$. Let $AK_{\omega}^G$ be the space of $\omega$-compatible $G$-invariant almost-complex structures. As $AK_{\omega}$ is contractible, it is also
connected; the same is true for $AK_{\omega}^G$ by taking the average of a path of metrics in $AK_{\omega}$ over $G.$

% and $J\in AK_{\omega}^T.$ 
%As in K\"ahler case (see \cite{fut-mab}), we obtain, using Moser Lemma, an analogue of the {\it extremal vector field} $Z^T_{\Omega}$ 

In this context, the following remark, generalizing Lemma $2$ in \cite{apo-cal-gau-Fri 2} suggests the definition of an {\it extremal
vector field} of $AK_{\omega}^G$ :
\begin{lemma}\label{lem1}
Let $J_t$ be a smooth family of almost-complex structures
compatible with the fixed symplectic form $\omega,$ which are
invariant under a compact group $G$ of symplectomorphisms
acting in a hamiltonian way on the compact symplectic manifold
$(M^{2n},\omega).$ %Denote by $\mathfrak{g}_{\omega}\in
%C^\infty(M)$ the finite dimensional vector space of smooth
%functions $f$ such that $X=grad_\omega f\in\mathfrak{g},$
%where $\mathfrak{g}$ denotes the Lie algebra of $G.$
Then, the
 $L^2$-orthogonal projection of the hermitian scalar curvature $s^{\nabla_t}$ of 
$(\omega,J_t,g_t)$ to $\mathfrak{g}_{\omega}$ is
independent of $t$.
\end{lemma}
\begin{proof}
By definition, any $f\in\mathfrak{g}_{\omega}$ defines a vector
field $X=grad_{\omega}f$ which is in $\mathfrak{g}$, and is therefore
Killing with respect to any of the metrics
$g_t(\cdot,\cdot)=\omega(\cdot,J_t\cdot)$ in $AK_{\omega}^G$. To prove our claim, we have to show
that $\int_M fs^{\nabla_t}\frac{\omega^n}{n!}$ is independent of $t$.
%Recall that, as a direct consequence of the Mohsen formula
%(see~\cite{gauduchon}), we have
%\begin{equation}\label{equa5}
%\frac{d}{dt}s^{\nabla^{g_t}}=-\delta J_t\delta \dot{J},\\
%\end{equation}
%where $\dot{J}=\frac{d}{dt}J_t$ and $\delta$ denote the
%codifferential with respect to $g_t$.
Using the relation (\ref{equa5}), we
obtain
\begin{equation*}
\frac{d}{dt}\left( \int_Mfs^{\nabla_t}\frac{\omega^n}{n!}
\right)=\int_M\left(-\delta^{g_t} J_t\delta^{g_t}
\dot{J}\right)f\frac{\omega^n}{n!}=\int_Mg_t(\dot{J},D^{g_t}d^c_{t}f)\frac{\omega^n}{n!}.
\end{equation*}
The fact that
$X=grad_\omega f$ is Killing implies that $D^{g_t}d^c_{t}f$ is an
anti-symmetric tensor. The result follows if we recall that $\dot{J}$
is a $g_t$-symmetric endomorphism of $T(M)$.
\end{proof}
\begin{remark}
The above lemma can also be viewed as a consequence from the fact that $s^{\nabla}$ is the momentum map for the action of $\textrm{Ham}(M,\omega)$. Indeed, consider a Lie subgroup $\mathbb{G}\subset\mathbb{H}$ where $\mathbb{H}$ is a Lie group (equipped with a bi-invariant metric) acting in a hamiltonian way on a symplectic manifold with moment map $\mu^\mathbb{H}$. Let $N$ be a $\mathbb{G}$-invariant connected subspace and denote by $\mu^\mathbb{G}$ the projection (with respect to the inner product) of the image of $\mu^\mathbb{H}$ on the dual of the Lie subalgebra of $\mathbb{G}$. Since $N$ is $\mathbb{G}$-invariant, the differential of $\mu^\mathbb{G}|_N$ (restriction of $\mu^\mathbb{G}$ to $N$) is zero. Therefore, $\mu^\mathbb{G}|_N$ is constant. In our case, $\mathbb{H}=\textrm{Ham}(M,\omega),\mathbb{G}=G, \textrm{Lie}(G)\cong\mathfrak{g}_{\omega}, N=AK_{\omega}^G$ and $\mu^\mathbb{H}$ is given by (\ref{equa35}). Lemma \ref{lem1} is equivalent to the fact that $\mu^\mathbb{G}|_N$ is constant.
\end{remark}
Given any $J\in AK_{\omega}^G$,  we define $z^G_{\omega}:=\Pi_{\omega} s^{\nabla},$ where $s^{\nabla}$ is the
hermitian scalar curvature of $(\omega,J,g)$. This $z^G_{\omega}$ is independant 
of $J$ by Lemma \ref{lem1}. Let $G=T$ be a maximal torus in $\textrm{Ham}(M,\omega)$. We obtain the following lemma
\begin{lemma}\label{lem10}
For any $J\in AK_{\omega}^T,$ the almost-K{\"a}hler metric $(\omega,J,g)$ is extremal if and only if 
\begin{equation*}
 \mathring{s}^{\nabla}=z^T_{\omega},
\end{equation*}
where $\mathring{s}^{\nabla}$ is the integral zero part of the hermitian scalar curvature ${s}^{\nabla}$ of $(\omega,J,g)$ given by $\mathring{s}^{\nabla}={s}^{\nabla}-\frac{\int_M{s}^{\nabla}\omega^n}{\int_M\omega^n}$.
\end{lemma}
\begin{proof}
Let $g$ be such an extremal metric, then $X=grad_{\omega}s^{\nabla}$ is a Killing field which is invariant by $T$. 
Denote by $\Xi=span \{X,\mathfrak{t}\} $ where $\mathfrak{t}=\textrm{Lie}(T)$, then the closure of $\{\exp\Xi\}$ in the (compact) 
isometry group is a compact torus which contains $T$. By the maximality of the torus, we have $X\in\mathfrak{t}$. The other direction is obvious.
\end{proof}
\begin{definition}\label{def1}
The vector field $Z^T_{\omega}:=grad_{\omega}z^T_{\omega}$ is called the {\it extremal vector field} relative to $T$. 
\end{definition}
\begin{remark}
The vector field $Z^T_{\omega}$ is invariant under isotopy of $\omega$ : let $\omega_t$ an isotopy of $T$-invariant symplectic forms in $\left[\omega\right]$ with $\omega_0=\omega,$ i.e.
\begin{equation*}
\omega_t=\omega_0+d\sigma_t,\;  0\leqslant t\leqslant1,
\end{equation*}
where $\sigma_t$ is a $T$-invariant $1$-form. The flow $\Phi_t$ of $X_t=-\sigma_t^{\sharp_{\omega_t}}$ ($\sharp_{\omega_t}$ stands for 
the isomorphism between $T^{\ast}(M)$ and $T(M)$ via $\omega_t$) verifies $(\Phi_t)^\ast\omega_t=\omega_0$ 
so $(\Phi_t)^\ast(Z^T_{\omega_t})=Z^T_{\omega_0}$. On the other hand, as the vector field $X_t$ is $T$-invariant, $(\Phi_t)^\ast(Z^T_{\omega_t})=Z^T_{\omega_t}$.

It follows that the introduced extremal vector field $Z^T_\omega$ coincides with the one defined by Futaki and Mabuchi in the K\"ahler case~\cite{fut-mab}.
Indeed, in a fixed  $T$-invariant K\"{a}hler class, any two K\"{a}hler forms  $\omega_1$ and $\omega_2$ are isotopic. Therefore, $Z^T_{\omega_1}=Z^T_{\omega_2}$.
\end{remark}

Denote by $\mathfrak{M}_{\omega}^T$ the set all $T$-invariant almost-K\"{a}hler metrics induced by $T$-invariant symplectic forms isotopic to $\omega.$ The extremal vector field $Z_{\omega}^T$ is an obstruction to the existence of metrics of constant hermitian scalar curvature in $\mathfrak{M}_{\omega}^T$. We recall that the space $\mathfrak{M}_{\omega}^T$ is related (via Moser's lemma) to the space $AK_{\omega}^T$.

\begin{corollary}\label{corol10}
If $\mathfrak{M}_{\omega}^T$ contains a metric with constant hermitian scalar curvature, then $Z_{\omega}^T=0$. Conversely, If $Z_{\omega}^T=0$, any extremal metric in $\mathfrak{M}_{\omega}^T$ is of constant hermitian scalar curvature.
\end{corollary}
\begin{proof}
We readily deduce from the defintion of $Z_{\omega}^T$ the first assertion of the corollary. Now, if $Z_{\omega}^T=0$, then $z^T_{\widetilde{\omega}}=0$ for any $\widetilde{\omega}$ isotopic to $\omega.$ In particular, any extremal metric in $\mathfrak{M}_{\omega}^T$ is of constant hermitian scalar curvature.
\end{proof}

\subsection{Periodicity of the extremal vector field}
In this section, we show the periodicity of the introduced vector field $Z_{\omega}^T$ relative to a fixed maximal torus $T$ in $\textrm{Ham}(M,\omega)$ when  $\left[\frac{\omega}{2\pi}\right]$ is integral modulo torsion.
\begin{theorem}
Assume that $\left[\frac{\omega}{2\pi}\right]$ is integral modulo torsion. Then, there exist a positive integer $\nu$ such that $\exp(2\pi{\nu} Z_{\omega}^T)=1.$
\end{theorem}
\begin{proof}
As in the previous section, we denote by
$\mathfrak{t}_{\omega}\subset C^\infty(M)$ the finite dimensional
space of smooth functions which are hamiltonians  with zero mean value of elements of
$\mathfrak{t}=\textrm{Lie}(T).$ %We fix an almost-complex structure $J\in AK_{\omega}^T$.

In the symplectic context, Futaki and Mabuchi defined in \cite{fut-mab1} a bilinear symmetric form $\Phi:\mathfrak{t}\times\mathfrak{t}\to \mathbb{R}$
\begin{equation}\label{equa14}
\Phi(X,Y)=\frac{1}{(2\pi)^2}\int_M {f^{X}}{f^{Y}}\frac{\omega^n}{(2\pi)^n},
\end{equation}
where $f^{X},f^{Y}\in\mathfrak{t}_{\omega}$ momentums of $X,Y$. %The extremal vector field defined in section \ref{sec3.2} is related to (\ref{equa14}) in the following way 

%By hypothesis, $\Omega$ is an integral cohomology class of $M$ modulo torsion. Hence, $(M^{2n},\omega,J)$ can be \textit{quantized} by a hermitian complex line bundle $(L,h)$.
%In other words, the curvature form $R^{\nabla^L}$, for some hermitian $\mathbb{C}$-linear connection $\nabla^L$ on $L$, is related to $\omega$ by $\frac{\sqrt{-1}}{2\pi}R^{\nabla^L}=\omega$.

By hypothesis, $\left[\frac{\omega}{2\pi}\right]$ is an integral cohomology class modulo torsion. Under this condition, Futaki and Mabuchi showed \cite{fut-mab1} that if $X,Y$ are generators of $S^1$-actions (i.e. $\exp\left(X\right)=\exp\left(Y\right)=1$), then
$\Phi(X,Y)\in\mathbb{Q}.$ To prove our claim, we have to show that
$2\pi\Phi(X_i,Z_{\omega}^T)\in\mathbb{Q}$ for all $1\leqslant i\leqslant k$, where $Z_{\omega}^T$ is the extremal vector field relative to $T$ and $X_1,\cdots,X_k$ are generators of the torus action. This would imply $2\pi Z_{\omega}^T\in\sum_{i=1}^k\mathbb{Q}X_i,$ so $\exp\left(2\pi{\nu}Z_{\omega}^T\right)=1$ for some positive integer $\nu$. 

%For any $X\in\mathfrak{t},$ we define  $\mathcal{F}_{\Omega}(X):=\Phi(X,Z^T_{\Omega})$
%=\int_Mf^X\tilde{s}^{\nabla^g}\omega^n$, where $f^X\in\mathfrak{t}_{\omega}$ is the momentum of $X$ and $\tilde{s}^{\nabla^g}$ is the integral zero part of ${s}^{\nabla^g}$ given by $\tilde{s}^{\nabla^g}={s}^{\nabla^g}-\frac{\int_M{s}^{\nabla^g}\omega^n}{\int_M\omega^n}.$

To show our claim, we recall Nakagawa's modified version~\cite{nagak} of Tian's formula~\cite{tian}
\begin{align}
\frac{{(n+1)!2^{n-1}}n!}{2\pi}\int_Mf\mathring{s}^{\nabla}\frac{\omega^n}{(2\pi)^n}&=\sum_{j=0}^n(-1)^j\left(\begin{array}{c}n \\j\end{array}\right)\int_M\Bigg{[}\left(\frac{1}{2}\frac{\Delta^g f}{2\pi}+\frac{\rho^{\nabla}}{2\pi}\right)+(N+n-2j)\left(\frac{f}{2\pi}+\frac{\omega}{2\pi}\right)\Bigg{]}^{n+1}\label{formula2}\\
&-\left(N+\frac{n\mu}{n+1}\right)2^n(n+1)!\int_M\left(\frac{f}{2\pi}+\frac{\omega}{2\pi}\right)^{n+1},\nonumber
\end{align}
%\begin{eqnarray}
%(n+1)!n!2^{n-1}\mathcal{F}_{\Omega}(X)&=&\sum_{j=0}^n(-1)^j\left(\begin{array}{c}n \\j\end{array}\right)\int_M\left[\left(-\frac{\Delta f^X}{2}+\rho^{\nabla}\right)+(n-2j)\left(f^X+\omega\right)\right]^{n+1}\nonumber\\
%&-&\frac{n\mu}{n+1}2^n(n+1)!\int_M\left(f+\omega\right)^{n+1},\label{formula1}
%\end{eqnarray}
where $f$ is a smooth function on $M$, $\mathring{s}^{\nabla}$ is the zero integral part of the hermitian scalar curvature ${s}^{\nabla}$ of $(\omega,J,g)$, $\rho^\nabla$ is the hermitian Ricci form, $\mu=\frac{\int_M\rho^{\nabla}\omega^{n-1}}{\int_M{\omega^n}},$ $N$ is an integer and $\Delta^g$ is the Riemannian Laplacian with respect to $g.$ The formula (\ref{formula2}) takes in account the normalization (\ref{equa21}) and is a direct consequence of the fact that
$\int_M\left(\Delta^gf\right)\omega^n=0$ and the identities
%$\left\{\begin{array}{ll} \sum_{j=0}^n(-1)^j\left(\begin{array}{c}n \\j\end{array}\right)\left(n-2j\right)^k=0\;{\text{ if }}k<n{\text{ or }}k=n+1. \\\sum_{j=0}^n(-1)^j\left(\begin{array}{c}n \\j\end{array}\right)\left(n-2j\right)^n=2^nn!.  \end{array}\right.$
%\begin{eqnarray*}
\begin{align*}
\sum_{j=0}^n(-1)^j\left(\begin{array}{c}n \\j\end{array}\right)\left(n-2j\right)^k&=0\;{\text{ if }}k<n{\text{ or }}k=n+1,\\
\sum_{j=0}^n(-1)^j\left(\begin{array}{c}n \\j\end{array}\right)\left(n-2j\right)^n&=2^nn!.
%\end{eqnarray*}
\end{align*}

As $\Phi$ is non-degenerate and $2\pi\Phi(X,Z_{\omega}^T)=\frac{1}{2\pi}\int_Mf^Xz_{\omega}^T\frac{\omega^n}{(2\pi)^n}=\frac{1}{2\pi}\int_Mf^X\mathring{s}^{\nabla}\frac{\omega^n}{(2\pi)^n}$, we thus reduced the problem to show that the right hand side of (\ref{formula2}) is rational when
$f=f^X\in\mathfrak{t}_{\omega}$ is the momentum with respect to $\omega$ of a hamiltonian vector field $X$ generating an $S^1$-action.
We have then essentially two integrals in the right hand side of (\ref{formula2}) which are $\int_Mf^X\omega^n$ and $\frac{1}{2\pi}\int_M\left(\frac{1}{2}{\Delta^g f^X}+(N+n-2j){f^X}\right){\widetilde{\omega}}^n$ with $\widetilde{\omega}=\frac{\rho^{\nabla}}{2\pi}+\left(N+n-2j\right)\frac{\omega}{2\pi}.$ 

In order to compute the latter integral, we will use the localization formula (\ref{equa13}). Since $M$ is compact, we consider an integer $N$ large enough such that $\widetilde{\omega}$ is symplectic.
%Let $X\in\mathfrak{t}$ generating a $S^1$-action with momentum $f^X$ with respect to $\omega$.
By Lemma \ref{lem4}, the vector field $X$ verifies
\begin{equation*}
2\pi\widetilde{\omega}(X,\cdot)=-d\left(\frac{1}{2}{\Delta^g f^X}+(N+n-2j){f^X}\right). 
\end{equation*}

On the other hand, Futaki and Mabuchi showed \cite{fut-mab1} that if $\int_Mf^X{\omega^n}=0$ then $f^X(p)\in2\pi\mathbb{Q}$ for any fixed point $p\in M$ of the $S^1$-action. Moreover, we have $\left(\Delta^gf^X\right)(p)=-tr\left(D^gdf^X\right)_p=tr\left(JD^gX\right)_p$. Note that $\left(JD^gX\right)_p$ is a symmetric hermitian endomorphism which is the generator of the induced linear $S^1$-action on $T_p(M)$ (see e.g. \cite{duf-sal}). This implies that
the trace $tr\left(JD^gX\right)_p\in 2\pi\mathbb{Z}$ and therefore ${\left(\Delta^gf^X\right)(p)}\in 2\pi\mathbb{Z}$. By (\ref{equa13}) and using the power series of the exponential function in $h$, we obtain 
\begin{equation*}
\frac{1}{2\pi}\int_M\left(\frac{1}{2}{\Delta^g f^X}+(N+n-2j){f^X}\right){\widetilde{\omega}}^n\in\mathbb{Q}. 
\end{equation*}
This concludes the proof.
\end{proof}

\subsection{The Hermitian-Einstein condition}
\begin{definition}
An almost-K\"ahler metric $(\omega,J,g)$ is called {\it Hermitian-Einstein} if the hermitian Ricci form $\rho^{\nabla}$ is a (constant) multiple of the symplectic form $\omega$, i.e.
\begin{equation*}
\rho^{\nabla}=\frac{s^{\nabla}}{2n}\omega,
\end{equation*}
so the hermitian scalar curvature $s^\nabla$ is constant.
\end{definition}

\begin{corollary}
If $c_1(M,\omega)$ is a multiple of $\left[\omega\right]$ and $Z_{\omega}^T=0$, then an extremal almost-K\"ahler metric $(J,g)$ in $\mathfrak{M}_{\omega}^T$ is Hermitian-Einstein if and only if $\rho^{\nabla}$ is ${J}$-invariant.
\end{corollary}
\begin{proof}
Let $(\widetilde{\omega},{J},{g})$ be an extremal almost-K\"ahler metric in $\mathfrak{M}_{\omega}^T.$ By Corollary \ref{corol10} and since $Z_{\omega}^T=0$, we deduce that the hermitian scalar curvature $s^{{\nabla}}$ of $(\widetilde{\omega},{J},{g})$ is constant. To prove our claim, it is enough to show that
$\rho^{{\nabla}}$ is co-closed (and therefore harmonic), i.e. $\delta^{{g}}\rho^{{\nabla}}=0$, where $\delta^{{g}}$ is the codifferential with respect to ${g}$.
Indeed, as $\rho^{{\nabla}}$ and $\widetilde{\omega}$ be then two harmonic forms representing the same cohomology class up to a multiple, Hodge theory implies
there are equal up to the same multiple.

Denote by ${d}^c$ the twisted exterior differential with respect to ${J}$ and $\Lambda_{{\widetilde{\omega}}}$ the contraction by $\widetilde{\omega}$. By hypothesis $\rho^{{\nabla}}$ is ${J}$-invariant, then ${{d}}^c\rho^{{{\nabla}}}=0.$ 
Now, using the relations (\ref{equa16}) and (\ref{equa24}), we have $\delta^{{{g}}}\rho^{{{\nabla}}}=\left[\Lambda_{{\widetilde{\omega}}},{{d}}^c\right]\rho^{{{\nabla}}}=-{{d}}^c\Lambda_{{\widetilde{\omega}}}\rho^{{{\nabla}}}=-\frac{1}{2}{{d}}^c{s}^{{\nabla}}=0$ since ${s}^{{\nabla}}$ is constant. Therefore, $\rho^{{\nabla}}$ is co-closed. The other direction is obvious.
\end{proof}

\section{Explicit formula of the hermitian scalar curvature} 
Let $(M^{2n},\omega,J,g)$ be an almost-K\"ahler
manifold. By Darboux theorem, there exist coordinates $\{z_i,t_i\}$ defined on an open set $U$ such that  the symplectic form $\omega$ has the form $\omega=\sum_{i=1}^ndz_i\wedge dt_i$ on $U$; $\{z_i,t_i\}$ are called the Darboux coordinates. 
In this section, we give an explicit formula of the hermitian scalar curvature ${s}^{\nabla}$ of $(\omega,J,g)$ in terms of the coordinates $\{z_i,t_i\}.$
On $U$, the metric $g$ is of the form,
\begin{equation*}
g=\sum_{i,j=1}^n G_{ij}(z,t) dz_i\otimes dz_j+ H_{ij}(z,t) dt_i\otimes dt_j + P_{ij}(z,t)dz_i\odot dt_j,
\end{equation*}
where $G=(G_{ij}),H=(H_{ij})$ are symmetric positive-definite matrix-valued functions which satisfy the compatibility conditions $GH-P^2=Id$ and $HP={}^t\! PH$ (where ${}^t\! P$ denote the transpose of $P=(P_{ij})$).
%{}^T \! 
We define a local section $\Phi$ of the anti-canonical bundle $K^{-1}_J(M)$ by
$\Phi:=(K_1^{\flat}-\sqrt{-1}JK_1^{\flat})\wedge\cdots\wedge(K_n^{\flat}-\sqrt{-1}JK_n^{\flat})$ where $K_i=\frac{\partial}{\partial t_i}$ and $\flat$ stands for the isomorphism between $T(M)$ and $T^\ast(M)$ via $g$.
Let $\phi$ and $\psi$ the real $n$-forms such that
$\Phi=\phi+\sqrt{-1}\psi$ (in fact $\psi$ and $\phi$ are related in the following way
$\psi(JX_1,\cdots,X_n)=\phi(X_1,\cdots,X_n)$). We still denote by $\nabla$ the hermitian connection induced on $K^{-1}_J(M)$ by the canonical hermtian connection. For any vector
field $X$, we have
\begin{equation*}
\nabla_X\phi=a(X)\phi+b(X)\psi
\end{equation*}
for some real 1-forms $a,b$ (we can deduce that
$a=\frac{1}{2}d\ln|\phi|^2$). Moreover, the hermitian Ricci form $\rho^{\nabla}$
is given by $\rho^{\nabla}=db$. Now, span $\{K_1,\cdots,K_n\}$ is Lagrangian. Hence,
$\Phi(K_1,\cdots,K_n)=\phi(K_1,\cdots,K_n)=\det g(K_i,K_j)=\det H$ which implies that
\begin{equation*}
(\nabla_X\Phi)(K_1,\cdots,K_n)=(a(X)-ib(X))\Phi(K_1,\cdots,K_n)=(a(X)-ib(X))
\det H.
\end{equation*}
Let $\beta(X):=(\nabla_X\Phi)(K_1,\cdots,K_n)/\det H$, then $\beta=trace\left(H^{-1}\circ B\right)$ where $B(X)=\left(g(\nabla_XK_i,K_j)+\sqrt{-1}g(\nabla_XK_i,JK_j)\right)=\left(2g(\nabla_{X}K_i^{1,0},K_j)\right).$
In particular,
\begin{equation}\label{equa30}
B(X^{0,1})=\left(2g(\nabla_{X^{0,1}}K_i^{1,0},K_j)\right).
%&=&\left(g(\nabla_{X^{0,1}}K_i,K_j)+\sqrt{-1}g(\nabla_{X^{0,1}}K_i,JK_j)\right)\nonumber
\end{equation}
Using Proposition \ref{prop5}, we can compute, via (\ref{equa30}), the imaginary part of $B$ and therefore $\rho^{\nabla}$ and $s^{\nabla}.$ If we denote by $H^{-1}=(H^{ij})$ the inverse of $H=(H_{ij})$ and $G_{ij,k}=\frac{\partial G_{ij}}{\partial z_k}, G_{ij}^{,k}=\frac{\partial G_{ij}}{\partial t_k}$ etc, we obtain (details of calculations are left)
\begin{align}
%A_{ij}&=&-\mathcal{I}m(B_{ij})\\
%&=&\frac{1}{2}\sum_{k,l}\left[-P_{kj,l}H_{il}+P_{kj}^{,l}P_{li}-P_{lj}H_{il,k} \right. \\
%&+&\left. H_{lj}P_{li,k}-P_{lj}P_{kl}^{,i}+H_{lj}G_{kl}^{,i}\right]dz_k\\
%&+&\left[-H_{kj,l}H_{il}+H_{kj}^{,l}P_{li}-  P_{lj}H_{il}^{,k}\right. \\
%&+&\left.H_{lj}P_{li}^{,k}-P_{lj}H_{kl}^{,i}+H_{lj}P_{lk}^{,i}\right]dt_k.\\\Bigg\{\bigg\{\Big\{\Big\{\big\{
%&=&\sum_{i,j}H^{ij}A_{ij}\\
%&=&\frac{1}{2}\sum_{k,i}\left(-P_{ki,i}+P_{ii,k}+G_{ki}^{,i}\right)dz_k\\
%&+&\frac{1}{2}\sum_{i,j,k,l}H^{ij}\left[P_{kj}^{,l}P_{li}-P_{lj}H_{il,k}-P_{lj}P_{kl}^{,i} \right]dz_k\\
%&+&\frac{1}{2}\sum_{k,i}\left(P_{ik}^{,i}+P_{ii}^{,k}-H_{ki,i}\right)dt_k\\
%&+&\frac{1}{2}\sum_{i,j,k,l}H^{ij}\left[H_{kj}^{,l}P_{li}-P_{lj}H_{il}^{,k}-P_{lj}H_{kl}^{,i} \right]dt_k\\
\rho^{\nabla}=&\frac{1}{2}\sum_{k,i,l}\left(P_{ki,i}^{,l}-G_{ki}^{,il}+P_{il,k}^{,i}-H_{li,ik}\right)dz_k\wedge dt_{l}\label{formule1}\\
+&\frac{1}{2}\sum_{i,j,k,l,r}\bigg{[}\left(H^{ij}P_{lj}H_{il,k}\right)^{,r}-\left(H^{ij}P_{lj}H_{il}^{,r}\right)_{,k}\bigg{]}dz_k\wedge dt_r\nonumber\\
+&\frac{1}{2}\sum_{k,i,l}\left(P_{ki,il}-P_{ii,kl}-G_{ki,l}^{,i}\right)dz_k\wedge dz_l+\frac{1}{2}\sum_{i,j,k,l,r}\left(H^{ij}P_{lj}H_{il,k}\right)_{,r}dz_k\wedge dz_r\nonumber\\
+&\frac{1}{2}\sum_{k,i,l}\left(P_{il}^{,ik}+P_{ii}^{,lk}-H_{li,i}^{,k}\right)dt_k\wedge dt_l-\frac{1}{2}\sum_{i,j,k,l,r}\left(H^{ij}P_{lj}H_{il}^{,r}\right)^{,k}dt_k\wedge dt_r.\nonumber\\
s^{\nabla}=&\sum_{i,j}\left(-G_{ij}^{,ij}-H_{ij,ij}+P_{ij,j}^{,i}+P_{ji,i}^{,j}\right)+\sum_{i,j,k,l}\left(H^{ij}P_{lj}H_{il,k}\right)^{,k}-\sum_{i,j,k,l}\left(H^{ij}P_{lj}H_{il}^{,k} \right)_{,k}.\label{formule2}
\end{align}
\subsection{The toric case}
A {\it toric} symplectic manifold $(M^{2n},\omega)$ is a symplectic manifold equipped with an effective hamiltonian action of an $n$-dimensional torus $T$. It is generated by a family of hamiltonian vector fileds $\{K_1,\cdots,K_n\}$  which are linearly independant on a dense open set $M^{\circ}$ and satisfy the condition $\omega(K_i,K_j)=0$ for all $i,j$. It is well-known that the symplectic form $\omega$ has the following local
expression $\omega=\sum_{i=1}^ndz_i\wedge dt_i,$ where $z_i$ is the
momentum coordinate and $t_i$ is a local coordinate such that
$K_i=\frac{\partial}{\partial t_i}.$ So, any $\omega$-compatible $T$-invariant almost-K\"ahler metric $g$ has the local expression
\begin{equation}\label{form50}
g=\sum_{i,j=1}^n\Big{(} G_{ij}(z) dz_i\otimes dz_j+ H_{ij}(z) dt_i\otimes dt_j + P_{ij}(z)dz_i\odot dt_j\Big{)}.
\end{equation}
An almost-K\"ahler manifold $(M^{2n},\omega,g)$ is called {\it locally toric} if in local coordinates the symplectic form is written as $\omega=\sum_{i=1}^ndz_i\wedge dt_i$ and the almost-K\"ahler metric $g$ has the form (\ref{form50}).

We deduce from (\ref{formule1}) and (\ref{formule2}) the expressions of the hermitian Ricci form $\rho^{\nabla}$ and the hermitian scalar curvature $s^{\nabla}$
\begin{align}
\rho^{\nabla}=&\frac{1}{2}\sum_{i,k,l}-H_{li,ik}dz_k\wedge dt_l+\frac{1}{2}\sum_{k,i,l}\left(P_{ki,il}-P_{ii,kl}\right)dz_k\wedge dz_l\label{formule3}\\
+&\frac{1}{2}\sum_{i,j,k,l,r}\left(H^{ij}P_{lj}H_{il,k}\right)_{,r}dz_k\wedge dz_r,\nonumber\\
s^{\nabla}=&-\sum_{i,j}^nH_{ij,ij},\label{formule4}
\end{align}
generalizing the formula of $s^{\nabla}$ given by Abreu \cite{abr} in the integrable case and rediscovering the expression found by Donaldson \cite{don1}.

Now, we suppose that the (lagrangian) $g$-orthogonal distribution to the $T$-orbits is involutive. This condition is automatically satisfied in the K\"ahler case. Then, using Frobenius' theorem, there exist local coordinates $\{t_i\}$ such that $\{dt_1,\cdots,dt_n\}$ spann the annihilator of the orthogonal distribution to the $T$-orbits and $\omega=\sum_{i=1}^ndz_i\wedge dt_i$. Any such $\omega$-compatible $T$-invariant almost-K\"ahler metric $(\omega,g)$ has the local expression
\begin{equation}\label{form1}
g=\sum_{i,j=1}^n\Big{(} G_{ij}(z) dz_i\otimes dz_j+ H_{ij}(z) dt_i\otimes dt_j\Big{)} {\text{  and  }}\omega=\sum_{i=1}^ndz_i\wedge dt_i.
\end{equation}
Therefore, the expression (\ref{formule3}) of $\rho^{\nabla}$ further simplifies
\begin{equation*}
\rho^{\nabla}=-\frac{1}{2}\sum_{i,k,l}H_{li,ik}dz_k\wedge dt_l.
\end{equation*} 
\begin{remark}
We can deform an $\omega$-compatible toric metric of the diagonal form (\ref{form1}) by considering $H^\epsilon=H+\epsilon U$ and $G^\epsilon=(H^\epsilon)^{-1}$ for some non-negative definite matrix-valued function $U=(U_{ij}),$ say with
compact support on an open set and $\epsilon$ small enough. The fact that the right hand side of (\ref{formule4}) is an under-determined linear differential operator implies that there are infinite dimensional families of $\omega$-compatible extremal almost-K\"ahler metrics around an extremal almost-K\"ahler metric of the form (\ref{form1}). In particular, infinite dimensional families of non-integrable extremal almost-K\"ahler metrics do exist around any extremal K\"ahler toric metric.
\end{remark}

\begin{lemma}\label{lem7}
Let $(M^4,\omega,J_0,g_0)$ be an almost-K\"ahler
manifold such that $(\omega,g_0)$ have the form (\ref{form1}) on some open set $V$ of $M$. Then, there exists an infinite dimensional family 
of almost-K\"ahler metrics $(\omega,J_\epsilon,g_\epsilon),$ defined for a
sufficiently small $\epsilon$, such that $\rho^{{\nabla}_\epsilon}=\rho^{\nabla_0}$ for all $\epsilon$; here $\nabla_\epsilon$ is the canonical hermitian connection corresponding to $J_\epsilon$.
Moreover, if $(\omega,J_0,g_0)$ is K\"ahler, we obtain an infinite dimensional family of non-integrable $\omega$-compatible almost-K\"ahler metrics .

\end{lemma}

\begin{proof}
let ${H}^\epsilon_{ij}(z)=H_{ij}(z)+\epsilon U_{ij}(z)$, where $\epsilon$ is a real, $U_{ij}= f_{ij}(z_1)h_{ij}(z_2)$ with $f_{ij}=f_{ji},h_{ij}=h_{ji}$. Now, the condition $\rho^{\nabla_\epsilon}-\rho^{\nabla_0}=0$ gives the following system of O.D.E.'s 
\begin{equation}
\sum_{k=1}^2\left(f_{ik}(z_1)h_{ik}(z_2)\right)_{,kj}=0,
\end{equation}
which reduces to the relations $f_{12}^{\prime}=\alpha f_{22}$, $f_{11}^{\prime}=\beta f_{12}$, $h_{12}^{\prime}=-\beta h_{11}$ and $h_{22}^{\prime}=-\alpha h_{12}$ for
some (real) constants $\alpha,\beta$ which we assume non-zero.
Choosing $f_{11}$ and $h_{22}$ arbitrary, the above relations determine the remaining functions $f_{12},f_{22},h_{11},h_{12}$. In particular, for $f_{11}$ and $h_{22}$ with compact support on $V$, $U_{ij}$  has a compact support on $V$. This ensures that for a sufficiently small $\epsilon$, ${H}^\epsilon=({H}^\epsilon_{ij})$ is positive-definite. Letting $G^\epsilon=(H^\epsilon)^{-1},$ we obtain the $\omega$-compatible metric

$g_\epsilon=\left\{\begin{array}{ll}\sum_{i,j=1}^2 G^\epsilon_{ij}(z) dz_i\otimes dz_j+ H^\epsilon_{ij}(z) dt_i\otimes dt_j  & {\text{  on $V,$}} \\g_0 & {\text{  elsewhere.}} \end{array}\right.$

One can check directly that for a generic choice of $f_{11}$ and $h_{22}$, the corresponding almost-complex structure $J_\epsilon$ is non-integrable.
\end{proof}

\begin{corollary}\label{prop7}
Let $(M^4,\omega,J,g)$ be a K\"ahler-Einstein (complex) surface which is locally toric. Then, $\omega$ admits an infinite dimensional family of non-integrable, $\omega$-compatible almost-complex structures inducing Hermitian-Einstein almost-K\"ahler metrics.
\end{corollary}
\begin{example}\label{exa1}
The corollary applies to the toric K\"ahler-Einstein surfaces $\mathbb{CP}^2$, $\mathbb{CP}^1\times\mathbb{CP}^1$ and $\mathbb{CP}^1\#3\overline{\mathbb{CP}^2}$, but also to the locally symmetric
Hermitian-Einstein spaces $\mathbb{CH}^2/\Gamma$ and $\left(\mathbb{CH}^1\times\mathbb{CH}^1\right)/\Gamma.$
\end{example}
\subsection{Extremal almost-K\"ahler metrics saturating LeBrun's estimates}
On an almost-K\"ahler manifold $(M^4,\omega,J,g)$ of dimension $4$, the bundle of $2$-forms decomposes as
\begin{equation*}
\wedge^2(M)=\wedge^+(M)\oplus\wedge^-(M),
\end{equation*}
where $\wedge^{\pm}(M)$ correspond to the eigenvalue $(\pm1)$ under the action of the (Riemannian) {\it Hodge operator} $\ast$ (see e.g. \cite{bes}). This decomposition is related to the splitting (\ref{split1}) as follows
\begin{equation*}
\wedge^+(M)=\mathbb{R}\,.\,\omega\oplus\wedge^{J,-}(M) {\text{ and }}\wedge^-(M)=\wedge_0^{J,+}(M).
\end{equation*}
The (Riemannian) curvature $R$, viewed as a (symmetric) linear map of $\wedge^2(M),$ decomposes as follows
\begin{equation*}
R=\left(\begin{array}{cc}W^++\frac{s}{12}Id|_{\wedge^+(M)} & \frac{1}{2}\widetilde{r_0}|_{\wedge^-(M)} \medskip\\   \frac{1}{2}\widetilde{r_0}|_{\wedge^+(M)}  & W^-+\frac{s}{12}Id|_{\wedge^-(M)} \end{array}\right)
\end{equation*}
where $W^{\pm}$ are symmetric trace-free endomorphism of $\wedge^{\pm}(M)$ respectively acting trivially on $\wedge^{\mp}(M)$ and $\widetilde{r_0}$ is  the (symmetric) operator defined on $\wedge^2(M)$ as $\widetilde{r_0}(X\wedge Y)=r_0(X,\cdot)\wedge Y+X\wedge r_0(Y,\cdot)$; here $r_0=r-\frac{s}{4}g$ is the trace-free part of the {\it Ricci tensor} $r$ which is the trace of the Riemannian curvature $R$ and $s$ is the {\it Riemannian scalar curvature} defined as the trace of the Ricci tensor
(for more details see \cite{apo-arm-dra}). 
The tensor $W^+$ (resp. $W^-$) is called the {\it selfdual Weyl tensor} (resp. {\it anti-selfdual Weyl tensor}). 

We recall now a (weak version) of LeBrun's result in \cite[Proposition 2.2]{leb}.
\begin{proposition} 
Let $(M^4,\omega)$ be a compact symplectic $4$-manifold and $g\in AK_\omega$ be a $\omega$-compatible almost-K\"ahler metric. Then
\begin{equation}\label{ineq}
V^{\frac{1}{3}}\left(\int_M\Bigg{|}\left( \frac{2}{3}s+2e\right)_{-}\Bigg{|}^{{3}}\frac{\omega^2}{2}\right)^{\frac{2}{3}}\ge 32\pi^2\left|c_1^+\right|_{L^2},
\end{equation}
where $V=\int_M\frac{\omega^2}{2}$ is the total volume of $(M^4,\omega),$ $e(x)$ is the lowest eignevalue of $W^+$ at $x,$ for any real-valued function $f$ on $M : f_{-}(x)=\min\left(f(x),0\right)$
and $c_1^+$ denotes the self dual part of the $g$-harmonic $2$-form representing $c_1(M,\omega).$

Moreover, the equality holds in (\ref{ineq}) if and only if $g$ is an extremal almost-K\"ahler metric with negative constant hermitian scalar curvature $s^\nabla$ and at each point $\omega$ is an eigenform of 
$W^+$ correponding to its lowest eigenvalue. 
\end{proposition}
Corollary \ref{prop7}, applied to the (complex) surfaces
$\mathbb{CH}^2/\Gamma$ and $\left(\mathbb{CH}^1\times\mathbb{CH}^1\right)/\Gamma,$ provides examples of non-integrable extremal almost-K\"ahler metrics saturating the inequality (\ref{ineq}). Indeed, for a Hermitian-Einstein almost-K\"ahler metric, $\rho^\nabla$
is $J$-invariant and therefore $\rho^\ast=R(\omega)$ is (see Sec. \ref{ricci}). Since  $R(\omega)=W^+(\omega)+\frac{s}{12}\omega$, we deduce that $\omega$ belongs to the eigenspace of $W^+$ if and only if $\rho^{\nabla}$ is $J$-invariant.
Recall that in the K\"ahler case the Riemannian curvature $R$ and $\widetilde{r_0}$ act trivially on $\wedge^{J,-}(M)$, and the selfdual Weyl tensor $W^+$ decomposes as
\begin{equation*}
W^+=\left(\begin{array}{ccc} \frac{s}{6} & 0 & 0 \\0 & -\frac{s}{12}  & 0 \\0 & 0 &-\frac{s}{12} \end{array}\right).
\end{equation*}
For the K\"ahler-Einstein metric on $\mathbb{CH}^2/\Gamma$ and $\left(\mathbb{CH}^1\times\mathbb{CH}^1\right)/\Gamma,$ $s=s^\nabla$ is negative and thus $\omega$ belongs to the lowest eigenspace of $W^+$. Then, by Corollary \ref{prop7}, we obtain an infinite dimensional family of non-integrable Hermitian-Einstein almost-K\"ahler metrics
whose the almost-K\"ahler form belongs to the lowest eigenspace of $W^+$ with non-positive constant hermitian scalar curvature $s^\nabla$. Other examples of non-integrable almost-K\"ahler metrics saturating the inequality (\ref{ineq}) appear in \cite{sun}.

%\section{Applications}

%The extremal vector field introduced in \ref{sec3.2} coincides with the one defined by Futaki and Mabuchi in the K\"ahler case~\cite{fut-mab}. A direct consequence of proposition \ref{prop1} is the following corollary
%\begin{corollary}\label{corol6}
%Suppose that $\omega_1$ and $\omega_2$ are two K\"ahler forms in the same K\"ahler class. Then $Z^T_{\omega_1}=Z^T_{\omega_2}$
%\end{corollary}
\section{further remarks and questions}

\;\;\;\;\;\;(1) In \cite{don1}, it is shown that in the toric case the existence of 
an extremal almost-K\"ahler metric is closely related (and 
conjecturally equivalent) to the existence of an extremal K\"ahler 
metric; this was further generalized to certain toric bundles in 
\cite{apo-cal-gau-Fri 2}. It will be interesting to establish a similar link in
general, assuming that there are integrable complex structures in $AK_{\omega}^T$.

(2) One expects that there will be infinite families of extremal almost-K\"ahler metrics (should they exist) in $AK^{T}_{\omega}$, thus 
generalizing the results of \cite{kim-sun} and Example \ref{exa1}. 
Indeed, this is suggested from the GIT formal picture in \cite{don}, where 
the existence and uniqueness (modulo the action of $\textrm{Ham}(M, \omega)$) of 
the extremal almost-K\"ahler metrics are expected to hold within a `stable' 
complexified orbit for the action of $\textrm{Ham}(M,\omega)$ on $AK_{\omega}$. 
However, as $\textrm{Ham}(M,\omega)$ does not have a nature complexification, the description of the complexified orbits is not obvious and 
is clearly established only on K\"ahler manifolds with $H^1(M,\mathbb{R})=\{0\}$ for the action of $\textrm{Ham}(M,\omega)$ on 
the subset $K_{\omega}$ of integrable, $\omega$-compatible almost-complex 
structures: In this case, the complexified orbit of $J \in K_{\omega}$ 
is identified in \cite{don} (via Moser's lemma) with the space of K\"ahler 
metrics within the K\"ahler class $\Omega=[\omega]$ on $(M,J)$. In the 
absence of holomorphic vector fields on $(M,J)$, the stability theorem of 
\cite{leb-sim} and \cite{fuj-sch} implies that the existence of an extremal K\"ahler metric in a given complexified orbit is an open condition on 
the space of such orbits; this was generalized in \cite{apo-cal-gau-Fri 2}, 
by fixing a maximal tourus $T \subset Aut_0(M,J) \cap \textrm{Ham}(M,\omega)$ 
(where $Aut_0(M,J)$ is the connected component of the finite dimensional 
group of complex automorphisms of $(M,J)$) and considering $T$-invariant
$\omega$-compatible complex structures. We expect the 
openess result to generally hold in the non-integrable case, by 
considering suitable complexified orbits in $AK_{\omega}^T$.

(3) A well-know result of Calabi states that any extremal K\"ahler 
metric $(\omega,g)$ on a compact complex manifold $(M,J)$ is invariant under a 
maximal connected compact subgroup $G$ of $Aut_0(M,J)\cap \textrm{Ham}(M,\omega)$. 
It will be desirable to know whether or not such a $G$ (or a maximal 
torus in it) is also maximal as a compact subgroup of $\textrm{Ham}(M,\omega)$? 
More generally, one would like to know whether or not an extremal 
almost-K\"ahler metric is necessarily invariant under a maximal torus in 
$\textrm{Ham}(M,\omega)$?

(4) It would be interesting to know how does the extremal vector field 
$Z_\omega^T$, which we introduced as an invariant of a maximal torus $T$
in $\textrm{Ham}(M,\omega)$, characterize this torus up to congugacy in $\textrm{Symp}(M, \omega)$ or $\textrm{Ham}(M,\omega)$; elaborating on the theory
appearing in recent works \cite{abr-gra-kit,kar-kes-mar} would be an 
appealing direction of further investigation.

\subsection*{Acknowledgements} The author thanks V. Apostolov for his judicious advice, and M. Pinsonnault for stimulating discussions.

%So the {\it Futaki character} is given by $\mathcal{F}_{\Omega}(X):=\Phi(X,Z^T_{\Omega})$. %So for $X\in\mathfrak{t}$, $\mathcal{F}_{\Omega}(X)=\int_Mf^X\mathring{s}^{\nabla}\frac{\omega^n}{n!}$ where $f^X$ is the momentum of $X$ with respect to $\omega$ and $\mathring{s}^{\nabla}=s^{\nabla}-\frac{\int_Ms^{\nabla}\frac{\omega^n}{n!}}{\int_M\frac{\omega^n}{n!}}$ the part of $s^{\nabla}$ wich integrates to zero. 

\end{document}